\pdfoutput=1
\documentclass[12pt,a4paper]{article}
\usepackage[margin=1in]{geometry}
\usepackage[T1]{fontenc}
\usepackage[utf8]{inputenc}
\usepackage{lmodern}
\usepackage{amsmath,amssymb,amsfonts,amsthm}
\usepackage{mathrsfs}
\usepackage{needspace,enumitem}
\usepackage{microtype}
\usepackage[numbers]{natbib}
\usepackage[hidelinks]{hyperref}
\numberwithin{equation}{section}
\theoremstyle{plain}
\newtheorem{theorem}{Theorem}[section]
\newtheorem{proposition}[theorem]{Proposition}
\newtheorem{lemma}[theorem]{Lemma}
\newtheorem{corollary}[theorem]{Corollary}
\theoremstyle{definition}
\newtheorem{example}[theorem]{Example}
\newtheorem{definition}[theorem]{Definition}
\theoremstyle{remark}

\setlist[enumerate,1]{label=(\arabic*),leftmargin=0pt,
 itemindent=2.2em,labelwidth=1.7em,labelsep=.5em,align=left,font=\normalfont,
 topsep=.18\baselineskip,itemsep=.14\baselineskip,parsep=0pt}
\setlist[enumerate,2]{label=(\alph*),leftmargin=0pt,
 itemindent=2.2em,labelwidth=1.7em,labelsep=.5em,align=left,font=\normalfont,
 topsep=.14\baselineskip,itemsep=.12\baselineskip,parsep=0pt}
\newcommand{\doi}[1]{\href{https://doi.org/#1}{\nolinkurl{https://doi.org/#1}}}
\hypersetup{pdftitle={Ideal Ball Topologies on Commutative Unital Quantales},
 pdfauthor={Zhang Yaning and Wei Yao},
 pdfsubject={Ideal ball topologies, continuity, and quantale quotients},
 pdfkeywords={quantale, lattice ideal, residual similarity, ideal ball topology, quantale quotient}}

\title{Ideal Ball Topologies on Commutative Unital Quantales}
\author{Zhang Yaning \quad Wei Yao\\[0.6em]
\small School of Mathematics and Statistics,\\
\small Nanjing University of Information Science and Technology,\\
\small Nanjing 210044, China\\[0.3em]
\small \href{mailto:yaowei0516@163.com}{\texttt{yaowei0516@163.com}} (Wei Yao)}
\date{}

\begin{document}
\raggedbottom
\frenchspacing
\maketitle

\begin{abstract}
Let $I$ be a lattice ideal of a commutative unital quantale $L$.
We prove that the balls
$C_I(x,r)=\{y\in L:s(x,y)\to r\in I\}$, where
$s(x,y)=(x\to y)\wedge(y\to x)$ and $r\in I$, form a base
for a topology on $L$. No involutivity or additional closure
condition on $I$ is required. Multiplication and residual are
separately continuous. For integral quantales, continuity of
multiplication at the unit, joint continuity, radius compatibility
and uniform ball refinement are equivalent. Under these
conditions, multiplication, residual, meet and join are uniformly
continuous. We also construct a Hausdorff ideal-ball topology on
a noninvolutive frame for which meet is not jointly continuous.
This topology cannot be represented by lattice-filter balls on
the same frame. Negation identifies ideal balls with filter balls
in the involutive integral case, and every principal ideal on an
integral quantale has a multiplicative-filter representation.
An ideal stabilizer describes indistinguishability. Under uniform
refinement, the quotient is a Hausdorff reflection and a
sobrification. It is a quantale quotient if and only if arbitrary
joins preserve indistinguishability. Principal ideals give
partition topologies, which are compact exactly when the quotient
is finite. We also state compactness and countability criteria
in terms of ideal balls.
\end{abstract}


\medskip

\section{Introduction}

Quantales were introduced by Mulvey \cite{Mulvey1986} and
provide an algebraic setting for residuation
\cite{GalatosJipsenKowalskiOno2007,Rosenthal1990}.
Their topological applications include metric spaces
\cite{Lawvere1973}, categorical metric structures
\cite{ClementinoTholen2003,Flagg1997}, approximation
\cite{HofmannWaszkiewicz2011}, and probabilistic metrics
\cite{HofmannReis2013}. Ordinary topologies induced by
quantale-valued metrics are studied in \cite{CookWeiss2021}.
In this paper, we study topologies on a commutative unital
quantale whose neighborhoods are defined by a lattice ideal
and the biresiduum. We study continuity of the operations,
separation properties
and quotient structures.

Ideal-induced topologies have been studied on MV-algebras and
residuated lattices. Hoo \cite{Hoo1997} constructs topological
MV-algebras from systems of MV-ideals and studies completion.
Wang and Zhao \cite{WangZhao2022} study ideal systems on
residuated lattices. In Section~5 of their paper, they also
construct distance-ball topologies from lattice ideals on
involutive residuated lattices. Their balls use
$(x\to y)*(y\to x)$ in place of our meet of the two residuals.
Proposition~\ref{prop:wang-comparison} shows that their topology
is finer in the common involutive integral setting. On complete
MV-algebras the balls coincide. Whether the inclusion can be
strict in general is left open here.
Nasr Esfahani, Khaksar Haghani and Heidarian
\cite{EsfahaniHaghaniHeidarian2022} define $I$-balls on
BL-algebras using BL-ideals. Here we require only a lattice
ideal and allow both noninvolutive and nonintegral quantales.

The construction is also related to filter topologies.
Luan and Yang \cite{LuanYang2017} introduced
Chang-distance filter topologies on MV-algebras. Subsequent work
concerns chains \cite{WuLuanYang2020}, topological
characterizations \cite{AsadzadehRezaeiJamalzadeh2020,LuanZhaoYang2021},
and joint continuity \cite{LuanWeberYang2021}.
For complete MV-algebras, our balls are
$\{y:d(x,y)\oplus r\in I\}$ and coincide with filter balls
under negation. Theorem~\ref{thm:involutive-bridge} proves this
correspondence for involutive integral quantales.
Proposition~\ref{prop:principal-filter-representation} gives a
filter representation for principal ideals on integral quantales,
without involutivity. For general ideals this may fail.
Theorem~\ref{ex:failure-new} gives a nonprincipal ideal on a
noninvolutive frame whose topology is Hausdorff but whose meet
is not jointly continuous. Every filter-ball topology on that
frame has jointly continuous meet, so no lattice filter gives
the same topology.

Our main continuity result is Theorem~\ref{thm:integral-new}.
For integral quantales, continuity of multiplication at the unit
is equivalent to joint continuity, radius compatibility and
uniform ball refinement. These conditions also imply uniform
continuity of multiplication, residual, meet and join.

Topological and uniform MV-algebras are studied further in
\cite{BarbieriWeber1998,Weber2012,NajafiRezaeiKouhestani2017,GanLuanDengYang2024,YinXieYang2025,LiYang2026}.
Related work on BL-algebras and residuated lattices concerns
linear topologies, continuity and completion
\cite{HeYangWang2021,RasouliDehghani2020,YangZhang2019,ZahiriBorzooei2016,YangXinHe2018,KhanegirRezaeiKouhestani2019,ChenXinYang2024,GeorgescuPopescu2005,XiaoYangXinMa2025,YangHeWang2026}.
For complete residuated lattices of finite order dimension,
Krupka \cite{Krupka2025} uses the biresiduum to define balls
with residuated Scott-open sets as radii. Here we use the radius sets
$D_r=\{a:a\to r\in I\}$. The two constructions therefore
use the same biresiduum with different families of radius sets.

Section~2 recalls the preliminaries. Section~3 constructs
the topology and compares it with filter-ball topologies.
Section~4 studies continuity. Section~5 treats separation and
quotients, including the Hausdorff reflection and sobrification
under uniform refinement
\cite{Isbell1964,GoubaultLarrecq2013,Johnstone1982,PicadoPultr2012}.
Here sobriety is understood in ordinary topology;
quantale-valued sobriety is studied in \cite{ZhangZhang2022}.
Section~6 gives compactness and countability criteria.

\section{Preliminaries}\label{sec:preliminaries}

We recall some algebraic and topological notions; see
\cite{GalatosJipsenKowalskiOno2007,Rosenthal1990,Isbell1964,Bourbaki1995,Engelking1989,Kelley1975}.

\subsection{Algebraic preliminaries}

\begin{definition}
An algebra $(L,\wedge,\vee,*,e,0,1)$ is called a
\emph{commutative unital quantale} if:
\begin{enumerate}
\item $(L,\wedge,\vee,0,1)$ is a complete lattice;
\item $(L,*,e)$ is a commutative monoid;
\item $*$ distributes over arbitrary joins.
\end{enumerate}
It is called \emph{integral} if $e=1$.
\end{definition}

Let $Q_1$ and $Q_2$ be commutative unital quantales, with units
$e_1,e_2$ and multiplications $*_1,*_2$, respectively.
A map $f:Q_1\to Q_2$ is a \emph{quantale homomorphism} if
\[
\begin{aligned}
f\left(\bigvee_{\lambda\in\Lambda}x_\lambda\right)
 &=\bigvee_{\lambda\in\Lambda}f(x_\lambda),\\
f(x*_1y)&=f(x)*_2f(y),\qquad f(e_1)=e_2
\end{aligned}
\]
for every indexed family $(x_\lambda)_{\lambda\in\Lambda}$ and
all $x,y\in Q_1$. The empty family is included, so $f$ preserves the least element.

An important special case of a quantale is a \emph{frame},
a complete lattice in which finite meets distribute over arbitrary
joins. With meet as multiplication, a frame is an integral
commutative quantale. Thus, for frames, the preceding notion of
unital quantale homomorphism is precisely that of a map preserving
arbitrary joins and finite meets.

Returning to a general commutative unital quantale, preservation
of arbitrary joins by multiplication gives a right adjoint
$a\to(-)$ to $a*(-)$, called the \emph{residual}. Thus
\[
a*c\leq b\Longleftrightarrow c\leq a\to b.
\]
With its residual, $L$ is a complete commutative residuated
lattice. Conversely, on a complete lattice the adjunction implies
preservation of arbitrary joins.

\begin{proposition}\label{prop:residuation-prelim}
Let $L$ be a commutative unital quantale. For $a,b,c\in L$,
the following properties hold.
\begin{enumerate}
\item Multiplication is monotone in both variables.
If $a\leq b$, then $b\to c\leq a\to c$ and
$c\to a\leq c\to b$.
\item $e\to a=a$, $a*(a\to b)\leq b$, and
$a\leq b$ if and only if $e\leq a\to b$.
\item The residual satisfies
\begin{equation}\label{eq:residual-composition}
(a*b)\to c=b\to(a\to c).
\end{equation}
\item If $L$ is integral, then $a*b\leq a\wedge b$.
For a frame with meet multiplication, the residual is
\[
a\to b=\bigvee\{c\in L:a\wedge c\leq b\}.
\]
\end{enumerate}
\end{proposition}

\begin{proof}
$(1)$ Suppose that $a\leq b$. Since $b=a\vee b$, join
preservation gives $b*c=(a*c)\vee(b*c)$, so $a*c\leq b*c$.
Commutativity gives monotonicity in the other variable.
Also $a*(b\to c)\leq b*(b\to c)\leq c$, whence
$b\to c\leq a\to c$. Since $c*(c\to a)\leq a\leq b$,
we also obtain $c\to a\leq c\to b$.

$(2)$ For every $d\in L$, residuation gives
$d\leq e\to a$ if and only if $e*d=d\leq a$; hence $e\to a=a$.
Applying the adjunction to $a\to b\leq a\to b$ yields
$a*(a\to b)\leq b$.
Finally, $e\leq a\to b$ is equivalent to $a*e=a\leq b$.

$(3)$ For each $d\in L$, associativity and residuation give
\[
\begin{aligned}
d\leq(a*b)\to c
&\Longleftrightarrow a*b*d\leq c\\
&\Longleftrightarrow b*d\leq a\to c\\
&\Longleftrightarrow d\leq b\to(a\to c).
\end{aligned}
\]
Thus the two residuals are equal.

$(4)$ If $e=1$, then $b\leq1$ gives $a*b\leq a*1=a$,
and similarly $a*b\leq b$. Hence $a*b\leq a\wedge b$.
For a frame, the right adjoint to $a\wedge(-)$ is the join of
all $c$ such that $a\wedge c\leq b$, which gives the stated formula.
\end{proof}

The residual also defines negation by $\neg a=a\to0$.
An integral quantale is called \emph{involutive} if
$\neg\neg a=a$ for all $a$.

We use the residual to define a similarity on $L$.

\begin{definition}\label{def:similarity-new}
Let $L$ be a commutative unital quantale. For $x,y\in L$, put
\[
s(x,y)=(x\to y)\wedge(y\to x).
\]
The element $s(x,y)$ is called the \emph{residual similarity},
or \emph{biresiduum}, of $x$ and $y$.
\end{definition}

By residuation, $e\leq s(x,y)$ if and only if $x=y$.
In particular, $s(x,x)\geq e$; equality with $1$ need not hold
when $L$ is nonintegral.

We first give some properties of the residual similarity.

\begin{proposition}\label{prop:sim-new}
Let $L$ be a commutative unital quantale. For $x,y,z,x',y'\in L$,
the following statements hold:
\begin{enumerate}
\item $s(x,y)=s(y,x)$;
\item $s(x,x)\geq e$;
\item $s(x,y)*s(y,z)\leq s(x,z)$;
\item $s(x,x')*s(y,y')\leq s(x*y,x'*y')$;
\item $s(x,x')*s(y,y')\leq s(x\to y,x'\to y')$;
\item $s(x,x')\wedge s(y,y')\leq s(x\wedge y,x'\wedge y')$;
\item $s(x,x')\wedge s(y,y')\leq s(x\vee y,x'\vee y')$.
\end{enumerate}
\end{proposition}

\begin{proof}
$(1)$ By the definition and commutativity of meet,
$s(x,y)=(x\to y)\wedge(y\to x)=s(y,x)$.

$(2)$ Since $x*e=x\leq x$, residuation gives $e\leq x\to x$.
Thus $s(x,x)=x\to x\geq e$.

$(3)$ Put $p=s(x,y)$ and $q=s(y,z)$. By residuation,
$x*p\leq y$, $y*p\leq x$, $y*q\leq z$ and $z*q\leq y$.
Consequently,
\[
x*p*q\leq y*q\leq z,\qquad
z*p*q\leq y*p\leq x.
\]
It follows that $p*q\leq(x\to z)\wedge(z\to x)=s(x,z)$.

$(4)$ Put $p=s(x,x')$ and $q=s(y,y')$.
Then $x*p\leq x'$, $x'*p\leq x$, $y*q\leq y'$ and $y'*q\leq y$.
Hence
\[
(x*y)*p*q=(x*p)*(y*q)\leq x'*y',\qquad
(x'*y')*p*q=(x'*p)*(y'*q)\leq x*y.
\]
Residuation gives $p*q\leq s(x*y,x'*y')$.

$(5)$ Put $p=s(x,x')$ and $q=s(y,y')$. We have
\[
x'*p*q*(x\to y)\leq x*q*(x\to y)\leq y*q\leq y',
\]
so $p*q*(x\to y)\leq x'\to y'$.
Similarly,
\[
x*p*q*(x'\to y')\leq x'*q*(x'\to y')\leq y'*q\leq y,
\]
so $p*q*(x'\to y')\leq x\to y$.
By residuation,
$p*q\leq s(x\to y,x'\to y')$.

$(6)$ Put $d=s(x,x')\wedge s(y,y')$.
Since $d\leq x\to x'$ and $d\leq y\to y'$, we have
$d*(x\wedge y)\leq d*x\leq x'$ and
$d*(x\wedge y)\leq d*y\leq y'$. Thus
$d*(x\wedge y)\leq x'\wedge y'$.
Using $d\leq x'\to x$ and $d\leq y'\to y$ in the same way gives
$d*(x'\wedge y')\leq x\wedge y$.
Residuation yields $d\leq s(x\wedge y,x'\wedge y')$.

$(7)$ Put $d=s(x,x')\wedge s(y,y')$.
Join preservation and residuation give
\[
d*(x\vee y)=(d*x)\vee(d*y)\leq x'\vee y',
\]
and
\[
d*(x'\vee y')=(d*x')\vee(d*y')\leq x\vee y.
\]
Hence $d\leq s(x\vee y,x'\vee y')$ by residuation.
\end{proof}

The following examples will be used below.

\begin{example}[Frames of open sets]\label{ex:frame-prelim}
For a topological space $X$, the set $\mathcal O(X)$ of open
subsets is a frame. Arbitrary joins are unions, finite meets
are intersections, and the meet of an arbitrary family is the
interior of its intersection. With $U*V=U\cap V$ and unit $X$,
its residual is
\[
U\to V=\operatorname{int}\bigl((X\setminus U)\cup V\bigr).
\]
When $X$ is discrete, $\mathcal O(X)=\mathcal P(X)$.
Thus the powerset Boolean algebra, with intersection as
multiplication and unit $X$, is an integral quantale.
\end{example}

\begin{example}[A nonintegral powerset quantale]
\label{ex:powerset-prelim}
Let $G$ be an abelian group, written additively.
The complete lattice $\mathcal P(G)$ is a commutative unital
quantale under
\[
A*B=\{a+b:a\in A,\ b\in B\},\qquad e=\{0\}.
\]
Its residual is $A\to B=\{g\in G:A+g\subseteq B\}$.
If $G$ has more than one element, its unit $\{0\}$ differs
from its greatest element $G$, so the quantale is nonintegral.
For instance,
\[
s(\{0\},\{0\})=\{0\}=e\ne G=1.
\]
Thus self-similarity need not equal the greatest element.
\end{example}

\begin{example}[Complete MV-algebras]\label{ex:MV-prelim}
Let $(A,\oplus,\neg,0)$ be a complete MV-algebra.
With
\[
1=\neg0,\qquad
x\odot y=\neg(\neg x\oplus\neg y),\qquad
x\to y=\neg x\oplus y,
\]
it is an involutive integral commutative quantale.
Writing $x\ominus y=x\odot\neg y$, its \emph{Chang distance} is
\[
d(x,y)=(x\ominus y)\oplus(y\ominus x).
\]
The residual similarity satisfies $s(x,y)=\neg d(x,y)$.
For example, the standard MV-operations on $[0,1]$ are
$x\oplus y=\min\{1,x+y\}$ and $\neg x=1-x$; hence
$x\odot y=\max\{0,x+y-1\}$ and $d(x,y)=|x-y|$.
Section~3 compares the corresponding ideal balls with Chang-distance balls.
\end{example}

We now recall the definition of a lattice ideal.
A \emph{lattice ideal} of $L$ is a nonempty downward closed
subset $I$ closed under binary joins. It is \emph{principal} if
$I=\mathord{\downarrow}c=\{r:r\leq c\}$, and
\emph{finitely generated} if it is the least lattice ideal
containing a finite nonempty subset. In particular, $0\in I$.
Throughout the paper, an ideal means a lattice ideal unless
otherwise stated. A \emph{lattice filter} is a nonempty upward
closed subset closed under binary meets; multiplicative closure
is not assumed.

We also recall quantale congruences and quantic nuclei.

\begin{definition}\label{def:quantale-congruence-nucleus-new}
Let $L$ be a quantale.
\begin{enumerate}
\item An equivalence relation $\theta$ on $L$ is a
\emph{quantale congruence} if it is compatible with multiplication
and arbitrary joins: $x\mathrel\theta x'$ and $y\mathrel\theta y'$
imply $x*y\mathrel\theta x'*y'$, and
\[
x_j\mathrel\theta y_j\ (j\in J)
\Longrightarrow
\bigvee_{j\in J}x_j\mathrel\theta\bigvee_{j\in J}y_j.
\]
\item A monotone, extensive and idempotent map $j:L\to L$ is
a \emph{quantic nucleus} if
$j(x)*j(y)\leq j(x*y)$ for all $x,y\in L$.
\end{enumerate}
\end{definition}

\subsection{Topological preliminaries}

We recall some notions from uniform spaces and separation theory.

\begin{definition}\label{def:uniform-prelim}
A \emph{uniformity} on a nonempty set $X$ is a proper filter
$\mathscr U$ on $X\times X$ such that, for every $U\in\mathscr U$,
\[
\Delta_X\subseteq U,\ U^{-1}\in\mathscr U,\
V\circ V\subseteq U\text{ for some }V\in\mathscr U.
\]
Here $\Delta_X$ is the diagonal, $U^{-1}$ is the inverse relation,
and $\circ$ denotes relational composition. The induced topology
has neighborhood base $\{U[x]:U\in\mathscr U\}$ at $x$, where
$U[x]=\{y:(x,y)\in U\}$.
\end{definition}
Cauchy nets, completeness and total boundedness
\cite{Isbell1964,Kelley1975} are expressed in terms of balls in
Definition~\ref{def:I-complete-tb}.

For the later discussion of sobriety, a nonempty closed subset
$C$ of a space is called \emph{irreducible} if
$C\subseteq C_1\cup C_2$, with $C_1,C_2$ closed, implies
$C\subseteq C_1$ or $C\subseteq C_2$. A \emph{generic point}
of $C$ is a point $x$ with $C=\overline{\{x\}}$.
A space is \emph{sober} if every nonempty irreducible closed set
has a unique generic point. Every Hausdorff space is sober and
every sober space is $T_0$ \cite{Johnstone1982,PicadoPultr2012}.
For the separation results below, we use complete regularity
without a separation axiom; a
completely regular Hausdorff space is called Tychonoff.
Finally, for the countability results, a space is \emph{first
countable}, denoted by $(A_1)$, if every
point has a countable neighborhood base. It is \emph{second
countable}, denoted by $(A_2)$, if its topology has a countable
base.

\section{Ideal-induced ball topologies}\label{sec:ideal-balls}

In this section, we define balls using the residual similarity
and a lattice ideal.

\begin{definition}\label{def:ball-new}
Let $L$ be a commutative unital quantale and let $I$ be a lattice
ideal of $L$. For $x\in L$ and $r\in I$, put
\[
C_I(x,r)=\{y\in L:s(x,y)\to r\in I\}.
\]
The set $C_I(x,r)$ is called the \emph{$I$-ball} with center
$x$ and radius $r$. The topology induced by the lattice ideal
$I$ is
\begin{equation}\label{eq:lattice-ideal-topology}
\tau_I=\{U\subseteq L:\forall x\in U\ \exists r\in I,
\ C_I(x,r)\subseteq U\}.
\end{equation}
By Theorem~\ref{thm:base-new}, $\tau_I$ is a topology with
the $I$-balls as a base.
\end{definition}

We first prove three properties of these balls.

\begin{lemma}\label{lem:refinement-new}
Let $I$ be a lattice ideal of a commutative unital quantale $L$.
Then the following statements hold:
\begin{enumerate}
\item $x\in C_I(x,r)$ for all $x\in L$ and $r\in I$;
\item If $r,t\in I$ and $r\leq t$, then $C_I(x,t)\subseteq C_I(x,r)$;
\item If $y\in C_I(x,r)$ and $t=s(x,y)\to r$, then
$C_I(y,t)\subseteq C_I(x,r)$.
\end{enumerate}
\end{lemma}

\begin{proof}
$(1)$ Since $s(x,x)\geq e$, we have
$s(x,x)\to r\leq e\to r=r\in I$.

$(2)$ The inequality $s(x,y)\to r\leq s(x,y)\to t$ and
downward closure of $I$ give the inclusion.

$(3)$ The assumption $y\in C_I(x,r)$ gives
$t=s(x,y)\to r\in I$, so $C_I(y,t)$ is defined.
If $z\in C_I(y,t)$, then $s(y,z)\to t\in I$.
By Proposition~\ref{prop:sim-new}(3) and
\eqref{eq:residual-composition},
\[
s(x,z)\to r
\leq\bigl(s(x,y)*s(y,z)\bigr)\to r
=s(y,z)\to t\in I.
\]
Thus $z\in C_I(x,r)$.
\end{proof}

The balls are decreasing with respect to the radius parameter:
larger radii give finer neighborhoods.

\begin{theorem}\label{thm:base-new}
Let $I$ be a lattice ideal of a commutative unital quantale $L$.
The family
\[
\mathscr C_I=\{C_I(x,r):x\in L,\ r\in I\}
\]
is a base for $\tau_I$, and $\{C_I(x,r):r\in I\}$ is a
neighborhood base at $x$.
\end{theorem}

\begin{proof}
By Lemma~\ref{lem:refinement-new}(1), every point belongs to a
ball. For the base intersection axiom, suppose that
$y\in C_I(x_1,r_1)\cap C_I(x_2,r_2)$. Then
$t_i=s(x_i,y)\to r_i\in I$ for $i=1,2$. Since $I$ is closed
under binary joins, $t=t_1\vee t_2$ also belongs to $I$.
Parts (2) and (3) of the lemma give
\[
C_I(y,t)\subseteq C_I(y,t_i)\subseteq C_I(x_i,r_i)
\qquad(i=1,2).
\]
Hence the balls form a base. By part (3), every basic open
set containing $y$ contains a ball centered at $y$.
These balls form a neighborhood base, and the induced topology
is the family in \eqref{eq:lattice-ideal-topology}.
\end{proof}

\subsection{Comparison with filter balls}

The filter construction replaces the residual membership test
by a product membership test. Let $F$ be a lattice filter on $L$
and define
\[
B_F(x,t)=\{y\in L:t*s(x,y)\in F\},
\qquad x\in L,\quad t\in F,
\]
and write $\tau_F$ for the topology generated by these balls.

These balls form a base. Indeed,
$t*s(x,x)\geq t*e=t$, so $x\in B_F(x,t)$.
If $y\in B_F(x,t)$, then $u=t*s(x,y)\in F$.
For $z\in B_F(y,u)$, the triangle inequality gives
\[
u*s(y,z)=t*s(x,y)*s(y,z)\leq t*s(x,z).
\]
Upward closure of $F$ therefore yields
$B_F(y,u)\subseteq B_F(x,t)$. For two balls containing the
same point, the meet of their refining radii belongs to $F$
and gives a ball contained in their intersection.

For involutive integral quantales, negation gives the following
relation between ideal balls and filter balls.

\begin{theorem}\label{thm:involutive-bridge}
Let $L$ be an involutive integral commutative quantale and let
$I$ be a lattice ideal. Then
$F=\neg I=\{\neg r:r\in I\}$ is a lattice filter, and
\[
C_I(x,r)=B_F(x,\neg r)\ (x\in L,\ r\in I).
\]
In particular, $\tau_I=\tau_F$.
\end{theorem}

\begin{proof}
Negation is an order-reversing involution, so it sends lattice
ideals to lattice filters. Residuation and involutivity give
\[
a\to r=\neg(a*\neg r).
\]
Since $u\in I$ if and only if $\neg u\in F$, we obtain
\[
s(x,y)\to r\in I
\Longleftrightarrow s(x,y)*\neg r\in F.
\]
Since every member of $F$ has the form $\neg r$, the result follows.
\end{proof}

On complete MV-algebras, the residual similarity is the negation
of the Chang distance, as shown below.

\begin{example}\label{ex:MV-specialization-new}
Let $A$ be a complete MV-algebra and let $I$ be a lattice ideal.
Since $s(x,y)=\neg d(x,y)$,
\[
C_I(x,r)=\{y:d(x,y)\oplus r\in I\}.
\]
Theorem~\ref{thm:involutive-bridge} identifies these balls with
the Chang-distance filter balls for $F=\neg I$
\cite{LuanYang2017,LuanWeberYang2021}.
They also agree with the distance balls in
\cite[Section 5]{WangZhao2022} on complete MV-algebras.
For general involutive residuated lattices, that paper uses
$d_W(x,y)=(x\to y)^*\oplus(y\to x)^*$, where
$u^*=u\to0$ and $u\oplus v=u^*\to v$.
Here $d_W(x,y)^*=(x\to y)*(y\to x)$ need not equal
$s(x,y)$. Example~\ref{ex:wang-similarity-distinction} shows that equality
can fail outside MV-algebras.
\end{example}

\begin{proposition}[Comparison with Wang--Zhao balls]\label{prop:wang-comparison}
Let $L$ be an involutive integral commutative quantale and let
$I$ be a lattice ideal. Write
\[
t(x,y)=(x\to y)*(y\to x),\qquad
W_I(x,r)=\{y:t(x,y)\to r\in I\}.
\]
These are the distance balls of \cite[Section 5]{WangZhao2022}.
For every $x\in L$ and $r\in I$,
\[
W_I(x,r)\subseteq C_I(x,r).
\]
Consequently, the Wang--Zhao topology is finer than $\tau_I$.
On complete MV-algebras the balls, and hence the topologies, coincide.
\end{proposition}

\begin{proof}
Involutivity gives $d_W(x,y)^*=t(x,y)$.
Integrality gives $t(x,y)\leq s(x,y)$, whence
$s(x,y)\to r\leq t(x,y)\to r$. Downward closure of $I$
proves the ball inclusion. Since both families are bases,
each $\tau_I$-open set is open in the Wang--Zhao topology.
On an MV-algebra both $t(x,y)$ and $s(x,y)$ equal the negation
of the Chang distance, proving equality.
\end{proof}

The next example shows that the two similarity formulas can
differ. It does not show that the topologies are different.

\begin{example}\label{ex:wang-similarity-distinction}
Let $R=\mathbb F_2[u,v]/(u^2,v^2)$ and let $L$ be its lattice
of ideals, with ideal multiplication. It is a finite integral
commutative quantale. Its seven ideals are
\[
0,\ Z=(uv),\ A=(u),\ B=(v),\ D=(u+v),\ M=(u,v),\ R.
\]
Negation is the annihilator: it interchanges $0,R$ and $Z,M$
and fixes $A,B,D$. Thus $L$ is involutive. The residual is
$J\to K=\{r\in R:rJ\subseteq K\}$, so
\[
A\to B=B\to A=M,\  s(A,B)=M,\ 
d_W(A,B)^*=M^2=Z<M.
\]
Thus $d_W(A,B)^*$ and $s(A,B)$ need not coincide in a finite
involutive quantale. This does not by itself imply that the
induced topologies differ.
\end{example}

For principal ideals, a filter representation exists under
integrality alone.

\begin{proposition}\label{prop:principal-filter-representation}
Let $L$ be an integral commutative quantale and let $c\in L$.
The set
\[
F_c=\{a\in L:a\to c\leq c\}
\]
is a multiplicatively closed lattice filter, and
\[
\tau_{\mathord{\downarrow}c}=\tau_{F_c}.
\]
No involutivity assumption is needed.
\end{proposition}

\begin{proof}
We have $1\in F_c$, and antitonicity of the residual makes
$F_c$ upward closed. If $a,b\in F_c$, then
\[
(a*b)\to c=b\to(a\to c)\leq b\to c\leq c.
\]
Thus $a*b\in F_c$. Since $a*b\leq a\wedge b$ in an integral
quantale, $F_c$ is also closed under binary meets.
By Lemma~\ref{lem:refinement-new}, $C_{\downarrow c}(x,c)$ is
the smallest ideal ball centered at $x$, and
\[
C_{\downarrow c}(x,c)=\{y:s(x,y)\in F_c\}.
\]
Fix $t\in F_c$. If $s(x,y)\in F_c$, multiplicative closure gives
$t*s(x,y)\in F_c$. Conversely, if $t*s(x,y)\in F_c$, then
$t*s(x,y)\leq s(x,y)$ by integrality, so upward closure gives
$s(x,y)\in F_c$. Hence
\[
B_{F_c}(x,t)=C_{\downarrow c}(x,c)\qquad(t\in F_c).
\]
Hence the two topologies coincide.
\end{proof}

The representing filter need not be the principal filter with
the same generator, as the following example shows.

\begin{example}\label{ex:ideal-filter-godel}
Let $L=[0,1]$ with $a*b=a\wedge b$, and fix $0<c<1$.
For $a\in[0,1]$, the residual is $a\to c=1$ when $a\leq c$
and $a\to c=c$ when $a>c$. Hence
$F_c=(c,1]$ in Proposition~\ref{prop:principal-filter-representation}.
The principal ideal $[0,c]$ therefore gives
\[
C_{[0,c]}(x,c)=\{y:s(x,y)>c\}.
\]
In contrast, the principal filter $[c,1]$ gives
\[
B_{[c,1]}(x,c)=\{y:s(x,y)\geq c\}.
\]
In the ideal topology, the smallest neighborhood of $c$ is
$\{c\}$. For the filter $[c,1]$, the smallest neighborhoods of
$c$ and $1$ are both $[c,1]$. Thus the principal ideal and principal filter generated by $c$
induce different topologies. The ideal topology is represented
by $F_c=(c,1]$.
\end{example}

\section{Continuity of quantale operations}\label{sec:continuity}

We first prove separate continuity, then give conditions for
joint continuity.

\subsection{Separate continuity}

\begin{theorem}\label{thm:separate-new}
Let $L$ be a commutative unital quantale, let $I$ be a lattice ideal
of $L$, and let $a\in L$. The maps
\[
x\longmapsto a*x,\  x\longmapsto a\to x,\ 
x\longmapsto x\to a
\]
are continuous on $(L,\tau_I)$. For each of these maps $f$,
\[
f(C_I(x,r))\subseteq C_I(f(x),r)\ (x\in L,\ r\in I).
\]
If $L$ is integral, the maps $x\mapsto a\wedge x$ and
$x\mapsto a\vee x$ have the same property.
\end{theorem}

\begin{proof}
Proposition~\ref{prop:sim-new} and $s(a,a)\geq e$ give
\[
s(x,y)\leq s(a*x,a*y),\ 
s(x,y)\leq s(a\to x,a\to y),\ 
s(x,y)\leq s(x\to a,y\to a).
\]
Thus, for each displayed map,
\[
s(f(x),f(y))\to r\leq s(x,y)\to r.
\]
For $y\in C_I(x,r)$, the right-hand side belongs to $I$,
and downward closure gives $f(y)\in C_I(f(x),r)$.
This proves the ball inclusion and hence continuity.
If $L$ is integral, $s(a,a)=1$; the meet and join estimates in
Proposition~\ref{prop:sim-new} give the same argument for the
two lattice maps.
\end{proof}

\subsection{Joint continuity}

To study joint continuity, put
\[
S(L)=\{s(x,y):x,y\in L\}
\]
for the set of similarity values.

\begin{definition}\label{def:RCs-new}\label{def:LCs-new}
\label{def:RC-new}\label{def:LC-new}
Let $L$ be a commutative unital quantale and let $I$ be a lattice
ideal of $L$. The pair $(L,I)$ is said to satisfy:
\begin{enumerate}
\item the \emph{radius compatibility condition} $(\mathrm{RC})$
if, for every $r\in I$, there exist $p,q\in I$ such that
\begin{equation}\label{eq:4-3}
a\to p,\ b\to q\in I\Longrightarrow (a*b)\to r\in I
\ (\forall a,b\in L);
\end{equation}
\item the \emph{lattice-radius compatibility condition}
$(\mathrm{LC})$ if, for every $r\in I$, there exist $p,q\in I$
such that
\begin{equation}\label{eq:4-4}
a\to p,\ b\to q\in I\Longrightarrow (a\wedge b)\to r\in I
\ (\forall a,b\in L).
\end{equation}
\end{enumerate}
Their \emph{similarity-level versions} $(\mathrm{RC}_s)$ and
$(\mathrm{LC}_s)$ are obtained by restricting $a,b$ to $S(L)$.
\end{definition}

Condition $(\mathrm{RC})$ will be used for multiplication and
residual, and $(\mathrm{LC})$ for meet and join.
Clearly, $(\mathrm{RC})$ implies $(\mathrm{RC}_s)$, and
$(\mathrm{LC})$ implies $(\mathrm{LC}_s)$.
If $L$ is integral, then $s(1,a)=a$, so $S(L)=L$.
Hence each condition is equivalent to its similarity-level version.

\begin{theorem}\label{thm:all-cont-new}
Let $L$ be a commutative unital quantale and let $I$ be a lattice
ideal of $L$. Then:
\begin{enumerate}
\item under $(\mathrm{RC}_s)$, multiplication and residual are
jointly continuous on $(L,\tau_I)$;
\item under $(\mathrm{LC}_s)$, meet and join are jointly continuous
on $(L,\tau_I)$.
\end{enumerate}
\end{theorem}

\begin{proof}
Let $(x,y)\in L^2$ and $r\in I$.

$(1)$ Choose $p,q\in I$ as in $(\mathrm{RC}_s)$.
If $u\in C_I(x,p)$ and $v\in C_I(y,q)$, then
$s(x,u)\to p,s(y,v)\to q\in I$, and therefore
\[
\bigl(s(x,u)*s(y,v)\bigr)\to r\in I.
\]
For $\diamond\in\{*,\to\}$, Proposition~\ref{prop:sim-new}
and antitonicity in the first residual argument give
\[
s(x\diamond y,u\diamond v)\to r
\leq \bigl(s(x,u)*s(y,v)\bigr)\to r.
\]
Downward closure of $I$ yields
$u\diamond v\in C_I(x\diamond y,r)$.

$(2)$ Choose $p,q\in I$ as in $(\mathrm{LC}_s)$.
If $u\in C_I(x,p)$ and $v\in C_I(y,q)$, then
\[
\bigl(s(x,u)\wedge s(y,v)\bigr)\to r\in I.
\]
For $\diamond\in\{\wedge,\vee\}$, the corresponding estimate in
Proposition~\ref{prop:sim-new} gives
\[
s(x\diamond y,u\diamond v)\to r
\leq \bigl(s(x,u)\wedge s(y,v)\bigr)\to r\in I.
\]
In each case, the required ball inclusion proves joint continuity.
\end{proof}

\begin{example}\label{ex:similarity-range-new}
Outside the integral case, $S(L)$ can be a proper subset of $L$.
Let $G=(\mathbb Z_2)^2$ and let $L=\mathcal P(G)$, with setwise
addition as multiplication and unit $\{0\}$. If $s(X,Y)$ is
nonempty, every $g\in s(X,Y)$ satisfies $X+g=Y$. Thus $s(X,Y)$
is a coset of the translation stabilizer of $X$ and has $1$, $2$
or $4$ elements. For distinct nonzero $u,v\in G$, the set
$T=\{0,u,v\}$ does not lie in $S(L)$.
Thus $(\mathrm{RC})$ includes a test with $a=b=T$ that is absent
from $(\mathrm{RC}_s)$. The difference in domains does not imply inequivalence of the conditions.
\end{example}

\subsection{Neighborhood criteria for joint continuity}

For $r\in I$, put
\begin{equation}\label{eq:4-5}
D_r=\{a\in L:a\to r\in I\}.
\end{equation}
Thus $C_I(x,r)=\{y:s(x,y)\in D_r\}$.
Each $D_r$ is an upper set containing $e$, and $r\leq t$ implies
$D_t\subseteq D_r$. For subsets $A,B\subseteq L$, write
$A*B=\{a*b:a\in A,b\in B\}$ and
$A\wedge B=\{a\wedge b:a\in A,b\in B\}$.

\begin{theorem}\label{thm:D-char-new}
Let $L$ be a commutative unital quantale and let $I$ be a lattice
ideal of $L$. Then:
\begin{enumerate}
\item $(\mathrm{RC}_s)$ holds if and only if, for every $r\in I$,
some $v\in I$ satisfies
\[
(D_v\cap S(L))*(D_v\cap S(L))\subseteq D_r;
\]
\item $(\mathrm{LC}_s)$ holds if and only if, for every $r\in I$,
some $v\in I$ satisfies
\[
(D_v\cap S(L))\wedge(D_v\cap S(L))\subseteq D_r;
\]
\item $(\mathrm{RC})$ holds if and only if, for every $r\in I$,
some $v\in I$ satisfies $D_v*D_v\subseteq D_r$;
\item $(\mathrm{LC})$ holds if and only if, for every $r\in I$,
some $v\in I$ satisfies $D_v\wedge D_v\subseteq D_r$.
\end{enumerate}
In all four conditions, the two input radii may be chosen equal.
\end{theorem}

\begin{proof}
$(1)$ Suppose that $(\mathrm{RC}_s)$ holds. Given $r\in I$,
choose $p,q\in I$ as in its definition and put $v=p\vee q$.
Since $D_v\subseteq D_p\cap D_q$, any $a,b\in D_v\cap S(L)$ satisfy
$a\to p,b\to q\in I$. Thus $(a* b)\to r\in I$, or
$a* b\in D_r$, proving the required inclusion.
Conversely, suppose that the displayed inclusion holds for some
$v\in I$ at each $r\in I$. If $a,b\in S(L)$ satisfy
$a\to v,b\to v\in I$, then $a,b\in D_v\cap S(L)$, so
$a* b\in D_r$. This is $(\mathrm{RC}_s)$ with $p=q=v$.

$(2)$ Suppose that $(\mathrm{LC}_s)$ holds. Given $r\in I$,
choose $p,q\in I$ as in its definition and put $v=p\vee q$.
Since $D_v\subseteq D_p\cap D_q$, any $a,b\in D_v\cap S(L)$ satisfy
$a\to p,b\to q\in I$. Thus $(a\wedge b)\to r\in I$, or
$a\wedge b\in D_r$, proving the required inclusion.
Conversely, suppose that the displayed inclusion holds for some
$v\in I$ at each $r\in I$. If $a,b\in S(L)$ satisfy
$a\to v,b\to v\in I$, then $a,b\in D_v\cap S(L)$, so
$a\wedge b\in D_r$. This is $(\mathrm{LC}_s)$ with $p=q=v$.

$(3)$ Suppose that $(\mathrm{RC})$ holds. Given $r\in I$,
choose $p,q\in I$ as in its definition and put $v=p\vee q$.
Since $D_v\subseteq D_p\cap D_q$, any $a,b\in D_v$ satisfy
$a\to p,b\to q\in I$. Thus $(a* b)\to r\in I$, or
$a* b\in D_r$, proving the required inclusion.
Conversely, suppose that the displayed inclusion holds for some
$v\in I$ at each $r\in I$. If $a,b\in L$ satisfy
$a\to v,b\to v\in I$, then $a,b\in D_v$, so
$a* b\in D_r$. This is $(\mathrm{RC})$ with $p=q=v$.

$(4)$ Suppose that $(\mathrm{LC})$ holds. Given $r\in I$,
choose $p,q\in I$ as in its definition and put $v=p\vee q$.
Since $D_v\subseteq D_p\cap D_q$, any $a,b\in D_v$ satisfy
$a\to p,b\to q\in I$. Thus $(a\wedge b)\to r\in I$, or
$a\wedge b\in D_r$, proving the required inclusion.
Conversely, suppose that the displayed inclusion holds for some
$v\in I$ at each $r\in I$. If $a,b\in L$ satisfy
$a\to v,b\to v\in I$, then $a,b\in D_v$, so
$a\wedge b\in D_r$. This is $(\mathrm{LC})$ with $p=q=v$.

In each case the two input radii may therefore be chosen equal.
\end{proof}

\subsection{Uniform refinement and continuity at the unit}

We now define uniform ball refinement and relate it to continuity.

\Needspace{10\baselineskip}
\begin{definition}\label{def:UR-new}
Let $L$ be a commutative unital quantale and let $I$ be a lattice
ideal of $L$. The pair $(L,I)$ is said to satisfy the
\emph{intrinsic uniform ball-refinement condition} $(\mathrm{UR})$
provided that, for every $r\in I$, there exists $v\in I$ such that
\begin{equation}\label{eq:uniform-refinement}
s(x,z)\to v,\ s(z,y)\to v\in I
\Longrightarrow s(x,y)\to r\in I\ (\forall x,y,z\in L).
\end{equation}
\end{definition}

By the triangle inequality, $(\mathrm{RC}_s)$ implies
$(\mathrm{UR})$.

\begin{proposition}\label{prop:RCs-UR-new}
Let $L$ be a commutative unital quantale and let $I$ be a lattice
ideal of $L$. If $(L,I)$ satisfies $(\mathrm{RC}_s)$, then it
satisfies $(\mathrm{UR})$.
\end{proposition}

\begin{proof}
For $r\in I$, choose $v\in I$ as in Theorem
\ref{thm:D-char-new}(1). If $s(x,z)\to v,s(z,y)\to v\in I$, then
$(s(x,z)*s(z,y))\to r\in I$. Since
$s(x,z)*s(z,y)\leq s(x,y)$, the antitonicity of the first argument
of the residual and the downward closure of $I$ give
$s(x,y)\to r\in I$.
\end{proof}

\begin{lemma}\label{lem:separation-uniformity-new}
Let $L$ be a commutative unital quantale and let $I$ be a lattice
ideal of $L$. Then the relations
\[
U_r=\{(x,y)\in L^2:s(x,y)\to r\in I\}\ (r\in I)
\]
form a base for a uniformity inducing $\tau_I$ if and only if
$(\mathrm{UR})$ holds. Denote this uniformity by $\mathscr U_I$.
If $(\mathrm{UR})$ holds and $s(x,y)\to r\notin I$, there is
$v\in I$ such that $C_I(x,v)\cap C_I(y,v)=\emptyset$.
\end{lemma}

\begin{proof}
Each $U_r$ is symmetric and contains the diagonal. Moreover,
$U_{r\vee t}\subseteq U_r\cap U_t$.
These relations form a uniformity base if and only if each $U_r$
contains $U_v\circ U_v$ for some $v\in I$
\cite{Isbell1964,Bourbaki1995,Kelley1975}; this is $(\mathrm{UR})$.
The induced topology is $\tau_I$, since $U_r[x]=C_I(x,r)$.
If $(x,y)\notin U_r$ and $U_v\circ U_v\subseteq U_r$, a point
in both $C_I(x,v)$ and $C_I(y,v)$ would contradict symmetry and
this inclusion.
\end{proof}

\Needspace{16\baselineskip}
For integral quantales, the converse also holds. We obtain
the following continuity theorem.

\begin{theorem}[Continuity and uniform refinement]\label{thm:integral-new}
Let $L$ be an integral commutative quantale and let $I$ be a
lattice ideal of $L$. The following conditions are equivalent:
\begin{enumerate}
\item multiplication is continuous at $(1,1)$;
\item multiplication is jointly continuous;
\item $(L,I)$ satisfies $(\mathrm{RC})$;
\item $(L,I)$ satisfies $(\mathrm{UR})$.
\end{enumerate}
Under these equivalent conditions, multiplication, residual, meet
and join are uniformly continuous for the ball uniformity
$\mathscr U_I$ and its product uniformity on $L^2$.
Independently, the following conditions are equivalent:
\begin{enumerate}
\item[(5)] meet is continuous at $(1,1)$;
\item[(6)] meet is jointly continuous;
\item[(7)] $(L,I)$ satisfies $(\mathrm{LC})$.
\end{enumerate}
Moreover, $(\mathrm{RC})$ implies $(\mathrm{LC})$.
\end{theorem}

\begin{proof}
$(1)\Rightarrow(3)$ Since $s(1,a)=a$, we have $C_I(1,r)=D_r$.
Continuity at $(1,1)$ therefore gives, for each $r\in I$,
$p,q\in I$ such that $D_p*D_q\subseteq D_r$.
This is exactly $(\mathrm{RC})$.

$(3)\Rightarrow(2)$ In the integral case, $S(L)=L$.
Thus $(\mathrm{RC})$ implies $(\mathrm{RC}_s)$, and
Theorem~\ref{thm:all-cont-new}(1) gives joint continuity.

$(2)\Rightarrow(1)$ Joint continuity implies continuity at $(1,1)$.

$(3)\Rightarrow(4)$ This follows from
Proposition~\ref{prop:RCs-UR-new}.

$(4)\Rightarrow(3)$ Fix $r\in I$ and choose $v\in I$
as in $(\mathrm{UR})$. Suppose that $a\to v,b\to v\in I$.
Integrality gives $a*b\leq a$, while residuation gives
$b\leq a\to(a*b)$. Hence $b\leq s(a,a*b)$, and
\[
s(1,a)\to v=a\to v\in I,\qquad
s(a,a*b)\to v\leq b\to v\in I.
\]
Applying $(\mathrm{UR})$ to $1,a,a*b$ yields $(a*b)\to r\in I$.
Thus $(\mathrm{RC})$ holds with both input radii equal to $v$.

$(5)\Rightarrow(7)$ Since $C_I(1,r)=D_r$, continuity of meet
at $(1,1)$ gives $D_p\wedge D_q\subseteq D_r$ for suitable
$p,q\in I$ at each $r\in I$. This is $(\mathrm{LC})$.

$(7)\Rightarrow(6)$ Condition $(\mathrm{LC})$ implies
$(\mathrm{LC}_s)$, so Theorem~\ref{thm:all-cont-new}(2) applies.

$(6)\Rightarrow(5)$ Joint continuity of meet implies continuity
at $(1,1)$.

To prove the final implication, suppose that $(\mathrm{RC})$
holds. For $r\in I$, choose $v\in I$ with $D_v*D_v\subseteq D_r$.
If $a,b\in D_v$, then $a*b\in D_r$ and $a*b\leq a\wedge b$.
Since $D_r$ is upward closed, $a\wedge b\in D_r$.
Thus $D_v\wedge D_v\subseteq D_r$, proving $(\mathrm{LC})$.

It remains to prove uniform continuity under (1)--(4).
By Lemma~\ref{lem:separation-uniformity-new}, $\mathscr U_I$
is a uniformity.
For $r\in I$, choose $p,q$ as in $(\mathrm{RC})$.
If $(x,x')\in U_p$ and $(y,y')\in U_q$, then
$\bigl(s(x,x')*s(y,y')\bigr)\to r\in I$.
Proposition~\ref{prop:sim-new}(4) therefore gives
$(x*y,x'*y')\in U_r$, proving uniform continuity of multiplication.
Proposition~\ref{prop:sim-new}(5) gives
$(x\to y,x'\to y')\in U_r$, proving uniform continuity of residual.
For meet, choose $p,q$ as in $(\mathrm{LC})$ and use
Proposition~\ref{prop:sim-new}(6) to obtain
$(x\wedge y,x'\wedge y')\in U_r$.
For join, the same input radii and Proposition~\ref{prop:sim-new}(7)
give $(x\vee y,x'\vee y')\in U_r$.
Thus all four operations are uniformly continuous.
\end{proof}

For a general, possibly nonintegral quantale, the same estimates
show that $(\mathrm{RC}_s)$ implies uniform continuity of
multiplication and residual, and that $(\mathrm{UR})$ together
with $(\mathrm{LC}_s)$ implies uniform continuity of meet and join.

We give two examples satisfying $(\mathrm{RC})$, followed by
a Hausdorff example in which meet is not jointly continuous.

\begin{example}\label{ex:boolean-ideal-balls-new}
Let $L$ be a complete Boolean algebra, with $a*b=a\wedge b$,
and let $I$ be a lattice ideal. If
$x\mathbin\triangle y=(x\wedge\neg y)\vee(y\wedge\neg x)$,
then
\[
C_I(x,r)=\{y:(x\mathbin\triangle y)\vee r\in I\}
=\{y:x\mathbin\triangle y\in I\}\ (r\in I).
\]
Thus the balls are the congruence classes modulo $I$ and do not
depend on the radius. Also
$D_r=\{a:\neg a\in I\}$, which is closed under meet.
Hence $(\mathrm{RC})$ and $(\mathrm{LC})$ hold, and all four
operations are jointly continuous.
\end{example}

\begin{example}\label{ex:product-chain-new}
Let $L=[0,1]$ have product multiplication, and let
$I=[0,c)$, where $0<c<1$. The residual is $a\to b=1$
when $a\leq b$, and $a\to b=b/a$ when $a>b$.
The residual similarity is
\[
s(x,y)=
\begin{cases}
1,&x=y=0,\\
\dfrac{x\wedge y}{x\vee y},&x,y>0,\\
0,&\text{otherwise}.
\end{cases}
\]
For every $r\in I$, $C_I(0,r)=\{0\}$. If $x>0$, then
\[
C_I(x,0)=(0,1],\ 
C_I(x,r)=\left(\frac r c x,\frac c r x\right)\cap(0,1]
\ (0<r<c).
\]
As $r$ tends to $c$ from below, these balls form a neighborhood
base for the usual topology on $(0,1]$. Thus $(L,\tau_I)$ is the topological sum of the
isolated point $0$ and the usual space $(0,1]$. It is Hausdorff
and sober but neither compact nor discrete.

For $r\in I$, choose $p=q$ with $\sqrt{rc}<p<c$. If
$a\to p,b\to p\in I$, the product residual gives
$a,b>p/c$, and hence $ab>p^2/c^2>r/c$.
It follows that $(ab)\to r\in I$, including when $r=0$.
Therefore $(\mathrm{RC})$ holds.
\end{example}

\subsection{A Hausdorff topology with no filter-ball representation}

The following ideal-ball topology has no filter-ball representation
on the same frame.

\begin{theorem}[Failure of filter-ball representation]\label{ex:failure-new}\label{ex:frame-ideal-failure}
There exist a noninvolutive frame $L$ and a nonprincipal lattice
ideal $I$ such that $(L,\tau_I)$ is Hausdorff, but meet is not
jointly continuous. Consequently, no lattice filter $F$ on $L$
satisfies $\tau_I=\tau_F$ for the filter balls of Section~3.1.
\end{theorem}

\begin{proof}
\emph{Construction.} Let $H=\{0,a,b,c,1\}$ be the five-element frame with
$0<a,b<c<1$, where $a$ and $b$ are incomparable,
$a\wedge b=0$ and $a\vee b=c$.
Its Heyting residual satisfies $a\to0=b$ and $b\to0=a$.
It is noninvolutive, since $\neg\neg c=1\ne c$.
Its residual table (row antecedent, column consequent) is
\[
\begin{array}{c|ccccc}
\to&0&a&b&c&1\\\hline
0&1&1&1&1&1\\
a&b&1&b&1&1\\
b&a&a&1&1&1\\
c&0&a&b&1&1\\
1&0&a&b&c&1
\end{array}
\]

On the frame $L=H^{\mathbb N}$, with pointwise multiplication
$x*y=x\wedge y$, consider the lattice ideal
\[
\begin{aligned}
I=\{x\in L:\ &x(k)\leq c\text{ for every }k,\\
&x(k)=0\text{ for all but finitely many }k\}.
\end{aligned}
\]
It is downward closed, and the join of two members has finite
support and is still bounded by $c$.
The ideal has no greatest element: given $p\in I$, choose a
coordinate outside its support and adjoin the value $c$ there.
Hence $I$ is nonprincipal.

\emph{Failure of joint continuity.} For arbitrary $p,q\in I$,
choose $k_0$ outside their supports.
Let $u,v\in L$ equal $1$ at every coordinate except $k_0$,
where $u(k_0)=a$ and $v(k_0)=b$.
For $k\ne k_0$, we have
$(u\to p)(k)=p(k)$ and $(v\to q)(k)=q(k)$.
At $k_0$, these values are respectively
$a\to0=b$ and $b\to0=a$.
Both residuals therefore have finite support and all their
coordinates are at most $c$, so they belong to $I$.
However,
\[
\bigl((u\wedge v)\to0\bigr)(k_0)=0\to0=1,
\]
and hence $(u\wedge v)\to0\notin I$.
Hence $(\mathrm{RC})$ fails for $r=0$.
Since multiplication is meet on this frame, $(\mathrm{LC})$
also fails. By Theorem~\ref{thm:integral-new}, meet is not
continuous at $(1,1)$, although it is separately continuous
by Theorem~\ref{thm:separate-new}.

\emph{Hausdorffness.} Although meet is not jointly continuous,
coordinate fibers separate points. For each $k$, let
$r^k\in I$ have value $c$ at $k$ and $0$ elsewhere.
In $H$, the inequality $t\to c\leq c$ holds exactly when $t=1$.
If $y\in C_I(x,r^k)$, the $k$th coordinate of
$s(x,y)\to r^k$ is at most $c$, so
$s(x(k),y(k))=1$. This gives $x(k)=y(k)$ and hence
\[
C_I(x,r^k)\subseteq\{y:y(k)=x(k)\}.
\]
This inclusion holds at every point of each coordinate fiber,
so the fibers are open. Distinct points are separated by fibers
at a coordinate where they differ.

\emph{Failure of filter representation.} Every filter-ball
topology on this frame has jointly continuous meet. Indeed,
for any lattice filter $F$ and $r\in F$, the set
$D^F_r=\{d:r\wedge d\in F\}$ is closed under meet.
If $u\in B_F(x,r)$ and $v\in B_F(y,r)$, meet closure and
the similarity estimate give
\[
r\wedge s(x\wedge y,u\wedge v)
\geq r\wedge s(x,u)\wedge s(y,v)\in F.
\]
This proves joint continuity of meet for filter balls, which
contradicts the failure established above if the two topologies coincide.
\end{proof}

\section{Separation and quotient spaces}\label{sec:separation}

Multiplication and residual preserve indistinguishability by
separate continuity. For meet and join, compatibility also holds
when $L$ is integral or $(\mathrm{LC}_s)$ holds. We study the
separation properties of the quotient and the conditions under
which it inherits a quantale structure.

\subsection{Indistinguishability}

We define the ideal stabilizer as follows.

\begin{definition}\label{def:indistinguishability-new}
Let $L$ be a commutative unital quantale and let $I$ be a lattice
ideal of $L$. Define the \emph{ideal stabilizer}
\[
K_I=\{a\in L:a\to r\in I\text{ for every }r\in I\}.
\]
For any $x,y\in L$, we say that $x$ and $y$ are
\emph{$I$-indistinguishable}, written $x\sim_Iy$, if
$s(x,y)\in K_I$. Equivalently, $y\in C_I(x,r)$ for every $r\in I$.
\end{definition}

\begin{lemma}\label{lem:KI-properties-new}
Let $L$ be a commutative unital quantale and let $I$ be a lattice
ideal of $L$. Then $K_I$ is an upward closed submonoid of $(L,*,e)$.
In particular, $1\in K_I$.
\end{lemma}

\begin{proof}
Since $e\to r=r\in I$ for every $r\in I$, we have $e\in K_I$.
If $a\in K_I$ and $a\leq b$, then $b\to r\leq a\to r\in I$, so
$b\in K_I$. In particular, $1\in K_I$. If $a,b\in K_I$, then
\[
(a*b)\to r=b\to(a\to r)\in I
\]
for every $r\in I$, since $a\to r\in I$. Thus $a*b\in K_I$.
\end{proof}

\begin{lemma}\label{lem:indistinguishable-new}
Let $L$ be a commutative unital quantale and let $I$ be a lattice
ideal of $L$. Then $x\sim_Iy$ if and only if $x$ and $y$
belong to exactly the same $\tau_I$-open sets. The relation $\sim_I$
is an equivalence relation, and every open set is saturated, that is,
a union of $\sim_I$-classes. Moreover,
\[
\overline{\{x\}}=[x]_I=\bigcap_{r\in I}C_I(x,r).
\]
Here $[x]_I$ denotes the $\sim_I$-class of $x$; we also write $[x]$
when the relation is clear. The quotient map
$q:L\longrightarrow L/{\sim_I}$, with the quotient topology on
its codomain, is open and continuous.
\end{lemma}

\begin{proof}
Suppose that $x\sim_Iy$ and that $O$ is open with $x\in O$.
Choose $r\in I$ with $C_I(x,r)\subseteq O$. Then $y\in O$.
Symmetry shows that $x$ and $y$ have the same open neighborhoods.
Conversely, if they have the same open neighborhoods, then
$y\in C_I(x,r)$ for every $r\in I$.
Thus $\sim_I$ is an equivalence relation and every open set is saturated.

If $y\in\overline{\{x\}}$, every $C_I(y,r)$ contains $x$.
By symmetry, $x\sim_Iy$. Conversely, if $x\sim_Iy$, every
neighborhood of $y$ contains $x$, so $y\in\overline{\{x\}}$.
For an open set $O$, saturation gives $q^{-1}(q(O))=O$.
Thus $q$ is open; continuity follows from the quotient topology.
\end{proof}

\begin{proposition}\label{prop:kolmogorov-reflection-new}
Let $L$ be a commutative unital quantale and let $I$ be a lattice
ideal of $L$. Then the quotient
\[
q:(L,\tau_I)\longrightarrow \overline L=L/{\sim_I}
\]
is the Kolmogorov reflection: $\overline L$ is $T_0$, and every
continuous map from $(L,\tau_I)$ to a $T_0$ space factors uniquely
through $q$. The inverse-image map is a frame isomorphism
\[
q^{-1}:\mathcal O(\overline L)\longrightarrow\mathcal O(L,\tau_I).
\]
\end{proposition}

\begin{proof}
If $[x]\ne[y]$, Lemma~\ref{lem:indistinguishable-new} provides an
open set containing one of $x,y$ but not the other. This set is
saturated, so its image under $q$ is open and distinguishes $[x]$
from $[y]$. Thus $\overline L$ is $T_0$.

Let $f:(L,\tau_I)\to X$ be continuous, where $X$ is $T_0$.
If $x\sim_Iy$, continuity implies that $f(x)$ and $f(y)$ belong
to the same open sets of $X$, and hence $f(x)=f(y)$. Therefore
$\overline f([x])=f(x)$ defines the unique map with
$f=\overline f\circ q$. It is continuous by the definition of
the quotient topology.

Inverse image preserves arbitrary unions and finite intersections.
Every open set $O\subseteq L$ is saturated, so
$q^{-1}(q(O))=O$, and $q(O)$ is open because $q$ is open.
Consequently $q^{-1}$ is a bijection between the two frames of
open sets and is a frame isomorphism.
\end{proof}

\subsection{Separation and topological quotients}

The next theorem gives separation criteria. The $T_0$ and $T_1$
criteria do not require $(\mathrm{UR})$.

\begin{theorem}\label{thm:separation-new}
Let $L$ be a commutative unital quantale and let $I$ be a lattice
ideal of $L$. The following statements are equivalent:
\begin{enumerate}
\item $(L,\tau_I)$ is $T_0$;
\item $(L,\tau_I)$ is $T_1$;
\item for all $x\ne y$ there is $r\in I$ with $s(x,y)\to r\notin I$;
\item $\sim_I$ is equality, equivalently
$\bigcap_{r\in I}C_I(x,r)=\{x\}$ for every $x\in L$.
\end{enumerate}
Under $(\mathrm{UR})$, $(L,\tau_I)$ is completely regular, and
these conditions are also equivalent to Hausdorffness and to the
Tychonoff property.
\end{theorem}

\begin{proof}
$(1)\Rightarrow(4)$ If $x\sim_Iy$, then Lemma
\ref{lem:indistinguishable-new} says that $x$ and $y$ belong to
the same open sets. Since the space is $T_0$, $x=y$.
The intersection formula in (4) follows from the same lemma.

$(4)\Rightarrow(3)$ If $x\ne y$, then $y$ does not belong to
$\bigcap_{r\in I}C_I(x,r)$. Thus some $r\in I$ satisfies
$s(x,y)\to r\notin I$.

$(3)\Rightarrow(2)$ Condition (3) implies $[x]_I=\{x\}$
for every $x$. By Lemma~\ref{lem:indistinguishable-new},
$\overline{\{x\}}=[x]_I$, so all singletons are closed.
Thus $(L,\tau_I)$ is $T_1$.

$(2)\Rightarrow(1)$ Every $T_1$ space is $T_0$.

Under $(\mathrm{UR})$, Lemma~\ref{lem:separation-uniformity-new}
gives a compatible uniformity. Hence the space is completely
regular, without a separation assumption \cite{Isbell1964,Kelley1975}.
If (3) holds, the same lemma provides disjoint balls around any
two distinct points, proving Hausdorffness. Conversely, every
Hausdorff space is $T_1$. Complete regularity then makes these
conditions equivalent to the Tychonoff property as well.
\end{proof}

We next give a criterion for sobriety.

\begin{corollary}\label{thm:sobriety-new}
Let $L$ be a commutative unital quantale and $I$ a lattice ideal.
The space $(L,\tau_I)$ is sober if and only if $\sim_I=\Delta_L$
and every nonempty irreducible closed subset is a singleton.
Under $(\mathrm{UR})$,
\[
\text{sober}\ \Longleftrightarrow\ T_0\ \Longleftrightarrow\ T_1
\ \Longleftrightarrow\ T_2\ \Longleftrightarrow\ \sim_I=\Delta_L.
\]
\end{corollary}

\begin{proof}
By Lemma~\ref{lem:indistinguishable-new}, a $T_0$ space of this
form has singleton point closures. Sobriety is therefore
equivalent to $T_0$ and the stated condition on irreducible closed
sets. Under $(\mathrm{UR})$, Theorem~\ref{thm:separation-new}
applies, and every Hausdorff space is sober.
\end{proof}

In the integral case, separation can be read directly from
the stabilizer.

\begin{corollary}\label{cor:integral-separation-new}
Let $L$ be an integral commutative quantale and let $I$ be a
lattice ideal of $L$. Then $(L,\tau_I)$ is $T_0$, equivalently
$T_1$, if and only if $K_I=\{1\}$.
\end{corollary}

\begin{proof}
In an integral quantale, $s(1,a)=a$, and $s(x,y)=1$ exactly when
$x=y$. Hence $K_I=\{1\}$ implies $T_1$ by Theorem
\ref{thm:separation-new}. Conversely, if $(L,\tau_I)$ is $T_0$
and $a\in K_I$, then $1\sim_Ia$, so $a=1$. Since $1\in K_I$,
we obtain $K_I=\{1\}$.
\end{proof}

\begin{example}\label{ex:separated-interval-new}
Let $L=[0,1]$ with the \L{}ukasiewicz t-norm
$a*b=\max\{0,a+b-1\}$ \cite{KlementMesiarPap2000}.
Fix $0<c\leq1$ and take the nonprincipal lattice ideal $I_c=[0,c)$.
Since $s(x,y)=1-|x-y|$ and $a\to r=\min\{1,1-a+r\}$,
\[
C_{I_c}(x,r)=\{y\in[0,1]:|x-y|<c-r\}\ (r\in I_c).
\]
Thus $\tau_{I_c}$ is the usual topology of $[0,1]$; it is compact,
Hausdorff, sober and nondiscrete. The radii
$r_n=c(1-1/n)$, $n\geq1$, give a countable neighborhood base.
For $a<1$, take $r=\max\{0,a+c-1\}\in I_c$. Then
$a\to r\geq c$, proving $K_{I_c}=\{1\}$.

For $r\in I_c$, choose $(r+c)/2<p<c$. If
$a\to p,b\to p\in I_c$, then $a,b>1+p-c$, and hence
$a*b>1+2p-2c>1+r-c$. Thus $(a*b)\to r<c$, proving
$(\mathrm{RC})$. The ideal $I_c$ is not closed under the
MV-operation $\oplus$: $(c/2)\oplus(c/2)=c\notin I_c$.
\end{example}

Separate continuity implies that multiplication and residual
preserve indistinguishability, as follows.

\begin{theorem}\label{thm:congruence-new}
Let $L$ be a commutative unital quantale and let $I$ be a lattice
ideal of $L$. Then $\sim_I$ is a congruence for multiplication and
residual. If $L$ is integral, or if $(\mathrm{LC}_s)$ holds,
then it is also a congruence for meet and join.
\end{theorem}

\begin{proof}
Continuous maps preserve topological indistinguishability.
Let $x\sim_Ix'$ and $y\sim_Iy'$. Applying the separately
continuous maps from Theorem~\ref{thm:separate-new} to one
variable at a time gives
\[
x\diamond y\sim_Ix'\diamond y\sim_Ix'\diamond y'
\qquad(\diamond\in\{*,\to\}).
\]
If $L$ is integral, apply the same argument to meet and join
using Theorem~\ref{thm:separate-new}. Under $(\mathrm{LC}_s)$,
the conclusion follows from joint continuity in
Theorem~\ref{thm:all-cont-new}.
\end{proof}

\begin{theorem}\label{thm:hausdorff-reflection-new}
Let $L$ be a commutative unital quantale and let $I$ be a lattice
ideal satisfying $(\mathrm{UR})$. Give
$\overline L=L/{\sim_I}$ the quotient topology.
\begin{enumerate}
\item $\overline L$ is Hausdorff and completely regular, and
$q:L\to\overline L$ is the Hausdorff reflection of $(L,\tau_I)$;
\item $q$ is a sobrification of $(L,\tau_I)$;
\item the formulas
\[
[x]*[y]=[x*y],\  [x]\to[y]=[x\to y],\  \overline e=[e]
\]
define separately continuous operations, and $q$ is an open
continuous unital monoid homomorphism;
\item if $(\mathrm{RC}_s)$ holds, multiplication and residual on
$\overline L$ are jointly continuous, and $q$ is the Hausdorff
reflection among topological commutative unital monoids.
If $(\mathrm{LC}_s)$ holds, the formulas
$[x]\wedge[y]=[x\wedge y]$ and $[x]\vee[y]=[x\vee y]$
also define jointly continuous operations.
\end{enumerate}
\end{theorem}

\begin{proof}
$(1)$ Distinct classes have disjoint saturated balls by Lemma
\ref{lem:separation-uniformity-new}; their open images separate
the classes. Thus $\overline L$ is Hausdorff.
If $D$ is closed in $\overline L$ and $[x]\notin D$,
complete regularity of $L$ gives a continuous map to $[0,1]$
separating $x$ from $q^{-1}(D)$. By Proposition
\ref{prop:kolmogorov-reflection-new}, this map factors through $q$,
and its factor separates $[x]$ from $D$.
Thus $\overline L$ is completely regular.
Every continuous map into a Hausdorff space factors uniquely
through $q$ by the same proposition, proving the reflection property.

$(2)$ By (1), $\overline L$ is Hausdorff and hence sober.
Every sober space is $T_0$, so Proposition
\ref{prop:kolmogorov-reflection-new} gives unique continuous
factorization through $q$ for every continuous map from $L$ to a
sober space. Thus $q$ is a sobrification.

$(3)$ Theorem~\ref{thm:congruence-new} makes multiplication and
residual well defined on the classes. Associativity, commutativity
and the unit law descend from $L$, and $q$ preserves the product
and unit. Lemma~\ref{lem:indistinguishable-new} makes $q$ open
and continuous. To prove separate continuity, fix $x$ and write
$m_x(y)=x*y$ and $\overline m_{[x]}([y])=[x*y]$.
The identity $\overline m_{[x]}\circ q=q\circ m_x$ and the
quotient property of $q$ show that $\overline m_{[x]}$ is continuous.
For residual, applying the same quotient argument to the continuous
maps $y\mapsto x\to y$ and $y\mapsto y\to x$ proves continuity
in each variable.

$(4)$ Under $(\mathrm{RC}_s)$, multiplication and residual on $L$
are jointly continuous by Theorem~\ref{thm:all-cont-new}(1).
Since $q$ is open and surjective, $q\times q$ is an open quotient
map. For either operation $\diamond$, the identity
\[
\overline\diamond\circ(q\times q)=q\circ\diamond
\]
therefore proves joint continuity on $\overline L$.
A continuous homomorphism into a Hausdorff topological monoid
factors uniquely through $q$ by (1); surjectivity of $q$ ensures
that its factor preserves multiplication and the unit.
Finally, under $(\mathrm{LC}_s)$, Theorem~\ref{thm:congruence-new}
makes meet and join well defined, and
Theorem~\ref{thm:all-cont-new}(2) makes them jointly continuous
on $L$. The same $q\times q$ argument proves joint continuity
of each induced lattice operation on $\overline L$.
\end{proof}

\subsection{Quantale quotients}

For the quotient to be a quantale quotient, arbitrary joins
must also preserve indistinguishability.

\begin{definition}\label{def:JC-new}
Let $L$ be a commutative unital quantale and let $I$ be a lattice
ideal. The pair $(L,I)$ is said to satisfy $(\mathrm{JC})$ if,
for every family $(x_\lambda,y_\lambda)_{\lambda\in\Lambda}$ in
$L\times L$,
\begin{equation}\label{eq:5-3}
\bigl(\forall\lambda\in\Lambda,\ x_\lambda\sim_Iy_\lambda\bigr)
\Longrightarrow
\bigvee_{\lambda\in\Lambda}x_\lambda
\sim_I\bigvee_{\lambda\in\Lambda}y_\lambda.
\end{equation}
\end{definition}

We call $(\mathrm{JC})$ \emph{join compatibility}.
It is defined without assuming $(\mathrm{UR})$.
Example~\ref{ex:finite-ideal-JC-new} shows that $(\mathrm{UR})$
does not imply $(\mathrm{JC})$.

\begin{theorem}\label{thm:quotient-quantale-new}
Let $L$ be a commutative unital quantale and let $I$ be a lattice
ideal. The following statements are equivalent:
\begin{enumerate}
\item $(\mathrm{JC})$ holds;
\item $\sim_I$ is a quantale congruence;
\item the formulas
\[
\bigvee_{\lambda\in\Lambda}[x_\lambda]
=\left[\bigvee_{\lambda\in\Lambda}x_\lambda\right],
\  [x]*[y]=[x*y]
\]
make $L/{\sim_I}$ a commutative unital quantale and make $q$ a
quantale homomorphism.
\end{enumerate}
Under these conditions, its residual is $[x]\to[y]=[x\to y]$.
\end{theorem}

The quotient joins determine the lattice order.
Arbitrary meets need not be represented by meets in $L$:
the quantale homomorphism $q$ preserves arbitrary joins,
multiplication and the unit, but need not preserve arbitrary meets.

\begin{proof}
$(1)\Rightarrow(2)$ Theorem~\ref{thm:congruence-new} gives
compatibility with multiplication, and $(\mathrm{JC})$ gives
compatibility with arbitrary joins. Thus $\sim_I$ is a quantale
congruence.

$(2)\Rightarrow(3)$ Under (2), the displayed joins and product
are well defined. The joins make the quotient a complete
join-semilattice, and hence a complete lattice, with order
\[
[a]\leq[b]\Longleftrightarrow[a\vee b]=[b].
\]
Associativity, commutativity, the unit law and distributivity over
arbitrary joins descend from $L$. The quotient map preserves
arbitrary joins, multiplication and the unit, as required.

$(3)\Rightarrow(1)$ If $x_\lambda\sim_Iy_\lambda$ for every
$\lambda$, then $q(x_\lambda)=q(y_\lambda)$. Join preservation gives
\[
q\left(\bigvee_\lambda x_\lambda\right)
=\bigvee_\lambda q(x_\lambda)
=\bigvee_\lambda q(y_\lambda)
=q\left(\bigvee_\lambda y_\lambda\right).
\]
This is $(\mathrm{JC})$. In particular, finite joins are well
defined without assuming $(\mathrm{LC}_s)$.

For the residual adjunction, if $[z]\leq[x\to y]$, then
$[x]*[z]\leq[x*(x\to y)]\leq[y]$. Conversely, suppose that
$[x]*[z]\leq[y]$, equivalently $[x*z\vee y]=[y]$.
Since $z\leq x\to(x*z\vee y)$ in $L$, residual compatibility gives
\[
[z]\leq[x\to(x*z\vee y)]=[x\to y].
\]
\end{proof}

We give a sufficient condition for $(\mathrm{JC})$ in terms
of $K_I$.

\begin{proposition}\label{prop:KI-meets-JC-new}
Let $L$ be a commutative unital quantale and let $I$ be a lattice
ideal. If $K_I$ is closed under arbitrary meets, then
$(\mathrm{JC})$ holds. This applies when $K_I=\mathord{\uparrow}k$
for some $k\in L$.
\end{proposition}

\begin{proof}
Suppose that $x_\lambda\sim_Iy_\lambda$ for every
$\lambda\in\Lambda$, and put
$d=\bigwedge_{\lambda\in\Lambda}s(x_\lambda,y_\lambda)\in K_I$.
Distributivity and residuation give
\[
\left(\bigvee_\lambda x_\lambda\right)*d\leq\bigvee_\lambda y_\lambda,
\ 
\left(\bigvee_\lambda y_\lambda\right)*d\leq\bigvee_\lambda x_\lambda.
\]
Thus $d\leq s(\bigvee_\lambda x_\lambda,\bigvee_\lambda y_\lambda)$.
The upward closure of $K_I$ proves $(\mathrm{JC})$.
\end{proof}

The quotient can also be described by a quantic nucleus.
For such a nucleus $j$, write $L_j=\{x\in L:j(x)=x\}$.
The fixed points form a quantale with the inherited order and
the following joins, multiplication and unit
\cite{Rosenthal1990,ZhangLaan2015}:
\[
\bigvee\nolimits_j x_\lambda=j\left(\bigvee_\lambda x_\lambda\right),
\  x*_j y=j(x*y),\  e_j=j(e).
\]

\begin{proposition}\label{prop:nucleus-quotient-new}
Let $L$ be a commutative unital quantale and let $I$ be a lattice
ideal satisfying $(\mathrm{JC})$. The join-preserving quotient map
$q$ has a right adjoint $q_*$. The composite $j_I=q_*q$ is a
quantic nucleus, and $L/{\sim_I}$ is isomorphic to its fixed-point
quantale. Conversely, a quantic nucleus whose fibers are the
$\sim_I$-classes implies $(\mathrm{JC})$.
\end{proposition}

\begin{proof}
Since $q$ preserves arbitrary joins, it has a right adjoint $q_*$.
The identity $qq_*=\mathrm{id}$ follows from surjectivity of $q$.
Thus $j_I=q_*q$ is monotone, extensive and
idempotent. Moreover,
\[
q\bigl(j_I(x)*j_I(y)\bigr)=q(x)*q(y)=q(x*y),
\]
whence $j_I(x)*j_I(y)\leq q_*q(x*y)=j_I(x*y)$ by adjunction.
The restrictions of $q$ and $q_*$ give inverse order isomorphisms
between $L_{j_I}$ and $L/{\sim_I}$ and preserve the indicated
operations. Conversely, $x\mapsto j(x)$ preserves arbitrary joins
into $L_j$. If its fibers are the $\sim_I$-classes, then
$x_\lambda\sim_Iy_\lambda$ for all $\lambda$ implies
$j(\bigvee_\lambda x_\lambda)=j(\bigvee_\lambda y_\lambda)$,
which is $(\mathrm{JC})$.
\end{proof}

For the principal ideal below, arbitrary joins preserve the
classes. By contrast, finite symmetric difference is preserved
by finitary Boolean operations but not by arbitrary unions.

\begin{example}\label{ex:godel-partition-new}
Let $L=[0,1]$ with $a*b=a\wedge b$, and fix $0\leq c<1$.
Take the principal ideal $I_c=[0,c]$. Here
\[
s(x,y)=\begin{cases}1,&x=y,\\x\wedge y,&x\ne y.\end{cases}
\]
Since $a\to r=1$ for $a\leq r$ and $a\to r=r$ for $a>r$,
$C_{I_c}(x,r)=\{y:s(x,y)>r\}$ for $r\leq c$. Its smallest ball is
\[
C_{I_c}(x,c)=
\begin{cases}\{x\},&x\leq c,\\(c,1],&x>c.\end{cases}
\]
Thus the topology is the partition topology with singleton blocks
in $[0,c]$ and one block $(c,1]$. The quotient is discrete and is
the sobrification. The original space is not $T_0$, unlike the
Hausdorff $t$-topology of \cite{Krupka2025}.
Taking both input radii equal to $c$ verifies $(\mathrm{RC})$,
so all finitary operations are continuous by Theorem
\ref{thm:integral-new}.
Arbitrary joins respect these classes. Indeed, if every member
of a family lies in $[0,c]$, equivalent representatives are equal
term by term. If some member exceeds $c$, the corresponding
representative also exceeds $c$, so both joins lie in the single
block $(c,1]$. Hence $(\mathrm{JC})$ holds, with nucleus
\[
j_c(x)=\begin{cases}x,&x\leq c,\\1,&x>c.\end{cases}
\]
\end{example}

\begin{example}\label{ex:finite-ideal-JC-new}
Let $L=\mathcal P(\mathbb N)$ be the Boolean quantale with product
given by intersection, and let $I$ be the ideal of finite subsets.
For $A,B\subseteq\mathbb N$ and finite $r$,
\[
s(A,B)\to r=(A\mathbin{\triangle}B)\cup r.
\]
Consequently every ball is
\[
C_I(A,r)=\{B\subseteq\mathbb N:A\mathbin{\triangle}B
\text{ is finite}\}.
\]
The topology is the partition topology modulo finite symmetric
difference, with discrete quotient. The finitary Boolean operations
preserve finite symmetric difference and are jointly continuous.
The condition $(\mathrm{JC})$ fails: for
$A_n=\{n\}$ and $B_n=\emptyset$ one has $A_n\sim_I B_n$ for every
$n$, but
\[
\bigcup_{n\in\mathbb N}A_n=\mathbb N\not\sim_I\emptyset
=\bigcup_{n\in\mathbb N}B_n.
\]
\end{example}

Without integrality, even binary joins can fail to preserve
the classes despite radius compatibility.

\begin{example}\label{ex:power-new}
Let $G=(\mathbb Z_2)^n$, $n\geq1$, and take $L=\mathcal P(G)$ with
setwise addition as multiplication. Its unit is $\{0\}$ and its
greatest element is $G$, so it is nonintegral. Put $I=\{\emptyset\}$.
The residual is
\[
X\to Y=\{a\in G:X+a\subseteq Y\}.
\]
Since $a+a=0$, one has $a\in s(X,Y)$ if and only if $Y=X+a$.
For $A\subseteq G$, $A\to\emptyset=\emptyset$ exactly when
$A\ne\emptyset$. Thus $C_I(X,\emptyset)$ is the translation orbit
of $X$, and the topology consists of unions of these orbits.

Since the only radius is $\emptyset$, take $p=q=\emptyset$ in
$(\mathrm{RC})$.
If $A\to\emptyset,B\to\emptyset\in I$, then $A$ and $B$ are
nonempty, as is their product. Hence $(A*B)\to\emptyset\in I$.
Thus $(\mathrm{RC})$ holds, and multiplication and residual are
jointly continuous.

Meet is not jointly continuous. For distinct $a,b\in G$, the
pairs $(\{a\},\{a\})$ and $(\{a\},\{b\})$ are topologically
indistinguishable in $L^2$, but their intersections $\{a\}$ and
$\emptyset$ belong to different orbits. Therefore
$(\mathrm{LC}_s)$, and hence $(\mathrm{LC})$, fails.
The original space is not $T_0$; its finite discrete orbit
quotient is its sobrification.

The condition $(\mathrm{JC})$ also fails, even for binary joins.
Choose $0\ne a\in G$ and put
$X_1=Y_1=X_2=\{0\}$ and $Y_2=\{a\}$.
Then $X_i\sim_IY_i$ for $i=1,2$, whereas $X_1\cup X_2=\{0\}$ and
$Y_1\cup Y_2=\{0,a\}$ are not translates. Hence
$(\mathrm{RC})$, and therefore $(\mathrm{UR})$, does not imply
$(\mathrm{JC})$.
\end{example}

\section{Compactness and countability}\label{sec:compactness}

In this section, we study compactness and countability.
Throughout, $L$ is a commutative unital quantale and $I$ is a
lattice ideal of $L$.

\subsection{Completeness, total boundedness and compactness}

We first define completeness and total boundedness using ideal balls.

\begin{definition}\label{def:I-complete-tb}
A net $(x_i)$ is \emph{$I$-Cauchy} if, for every $r\in I$,
$s(x_i,x_j)\to r\in I$ for all sufficiently large $i,j$.
The pair $(L,I)$ is \emph{$I$-complete} if every such net
converges in $\tau_I$, and \emph{$I$-totally bounded} if, for
each $r\in I$, some finite $K_r\subseteq L$ satisfies
\begin{equation}\label{eq:6-1}
L=\bigcup_{a\in K_r}C_I(a,r).
\end{equation}
\end{definition}

Under $(\mathrm{UR})$, the equality $U_r[x]=C_I(x,r)$ identifies
these notions with completeness and total boundedness for
$\mathscr U_I$. The standard uniform-space compactness theorem
\cite{Isbell1964} therefore gives the following result.

\begin{theorem}\label{thm:compactness-new}
Under $(\mathrm{UR})$, $(L,\tau_I)$ is compact if and only if
$(L,I)$ is $I$-complete and $I$-totally bounded. Under the
separation conditions of Theorem~\ref{thm:separation-new}, these
conditions are equivalent to compact Hausdorffness.
\end{theorem}

\begin{proof}
Under $(\mathrm{UR})$, Lemma~\ref{lem:separation-uniformity-new}
identifies $\tau_I$ with the topology of $\mathscr U_I$.
The equivalence follows from Definition~\ref{def:I-complete-tb}
and the compactness theorem for uniform spaces, which does not
require separation \cite{Isbell1964}. Under the separation conditions, the topology
is Hausdorff by Theorem~\ref{thm:separation-new}.
\end{proof}

Compactness passes to and from the indistinguishability quotient.
This does not require $(\mathrm{UR})$.

\begin{proposition}\label{prop:compact-reflection-new}
The space $(L,\tau_I)$ is compact if and only if $L/{\sim_I}$
is compact.
\end{proposition}

\begin{proof}
If $L$ is compact, its continuous image $L/{\sim_I}$ is compact.
Conversely, let $(O_\lambda)$ be an open cover of $L$.
Each $O_\lambda$ is saturated and $q$ is open, so
$(q(O_\lambda))$ is an open cover of the quotient.
A finite subcover of the quotient pulls back to a finite subcover
of $L$, because $q^{-1}(q(O_\lambda))=O_\lambda$.
\end{proof}

\subsection{Principal and finitely generated ideals}

For principal ideals, the topology is a partition topology.
This gives a simple compactness criterion.

\Needspace{13\baselineskip}
\begin{theorem}\label{thm:principal-compact-new}
Let $L$ be a commutative unital quantale and
$I=\mathord{\downarrow}c$. Define
\begin{equation}\label{eq:6-2}
x\equiv_c y\Longleftrightarrow s(x,y)\to c\leq c.
\end{equation}
Then $(L,I)$ satisfies $(\mathrm{RC})$, and the following hold:
\begin{enumerate}
\item $\equiv_c$ coincides with $\sim_I$, and
$C_I(x,c)=[x]_{\equiv_c}$ is the smallest neighborhood of $x$.
The open sets are exactly the unions of these equivalence classes;
\item $(L,I)$ is $I$-complete, and compactness and $I$-total
boundedness are each equivalent to finiteness of $L/{\equiv_c}$;
\item the quotient is discrete and $L\to L/{\equiv_c}$ is a
sobrification. The space $(L,\tau_I)$ is sober if and only if
$\equiv_c=\Delta_L$, equivalently if and only if $\tau_I$ is discrete.
\end{enumerate}
\end{theorem}

\begin{proof}
To verify $(\mathrm{RC})$, take $p=q=c$.
If $a\to c\leq c$, $b\to c\leq c$ and $r\in I$, then
\[
(a*b)\to r\leq(a*b)\to c=b\to(a\to c)\leq b\to c\leq c.
\]

$(1)$ The radius $c$ gives the smallest ball, since $r\leq c$
implies $C_I(x,c)\subseteq C_I(x,r)$. Hence
$U_c=\bigcap_{r\in I}U_r$, and Lemma
\ref{lem:indistinguishable-new} identifies it with $\sim_I$.
Every neighborhood of $x$ contains $C_I(x,c)$, so the topology
consists of unions of these classes.
Conversely, each class $C_I(x,c)$ is open, so every union of
classes is open.

$(2)$ Let $(x_i)$ be $I$-Cauchy. Applying the Cauchy condition
at radius $c$ shows that all sufficiently late terms lie in one
$\equiv_c$-class. Every neighborhood of a point in this class
contains the whole class by (1), so the net converges to every
point of the class. Hence $(L,I)$ is $I$-complete.

At radius $c$, a finite ball cover is exactly a finite cover by
equivalence classes. Thus $I$-total boundedness implies that there
are finitely many classes. Conversely, if there are finitely many
classes, their $c$-balls form a finite cover and refine the balls
at every other radius. The quotient is discrete by (1), so
Proposition~\ref{prop:compact-reflection-new} gives the same
finite-class criterion for compactness.

$(3)$ The discrete quotient is a sobrification by Theorem
\ref{thm:hausdorff-reflection-new}. A partition topology is $T_0$
exactly when all its classes are singletons, and
$\text{discrete}\Rightarrow\text{sober}\Rightarrow T_0$.
\end{proof}

For a powerset Boolean quantale, two sets lie in the same class
precisely when they agree outside the generator of the ideal.

\begin{example}\label{ex:boolean-trace-new}
Let $X$ be a set, and regard $L=\mathcal P(X)$ as the Boolean
quantale with multiplication given by intersection. Fix $C\subseteq X$
and take $I=\mathord{\downarrow}C$. Since
\[
A\to B=(X\setminus A)\cup B,\ 
s(A,B)=X\setminus(A\mathbin{\triangle}B),
\]
one obtains
\[
A\equiv_C B
\Longleftrightarrow A\mathbin{\triangle}B\subseteq C
\Longleftrightarrow A\cap(X\setminus C)=B\cap(X\setminus C).
\]
Thus two sets belong to the same smallest ball if and only if
their intersections with $X\setminus C$ coincide. The quotient map
\[
q_C:\mathcal P(X)\longrightarrow\mathcal P(X\setminus C),
\  q_C(A)=A\cap(X\setminus C)
\]
preserves arbitrary joins and multiplication. Hence $(\mathrm{JC})$
holds, and the quotient is the Boolean quantale
$\mathcal P(X\setminus C)$. Its associated nucleus is
$j_C(A)=A\cup C$.

By Theorem~\ref{thm:principal-compact-new}, $(L,\tau_I)$ is compact
if and only if $X\setminus C$ is finite, and sober if and only if
$C=\emptyset$. For $C=X$ the topology is indiscrete; for
$C=\emptyset$ it is discrete.
\end{example}

Since every finitely generated lattice ideal is principal,
we obtain the following corollary.

\begin{corollary}\label{cor:finitely-generated-compact-new}
Let $L$ be an infinite commutative unital quantale and let $I$ be a
lattice ideal. If $(L,\tau_I)$ is compact and Hausdorff, then
$I$ is not finitely generated.
\end{corollary}

\begin{proof}
If $I$ were generated by $\{c_1,\ldots,c_n\}$, then
$I=\mathord{\downarrow}(c_1\vee\cdots\vee c_n)$.
By Theorem~\ref{thm:principal-compact-new}, Hausdorffness would
make $\tau_I$ discrete, and compactness would force $L$ to be finite.
\end{proof}

\subsection{Countability}

We now consider first countability, second countability and
metrizability.

\begin{definition}\label{def:countable-radius-new}
Let $L$ be a commutative unital quantale and let $I$ be a lattice
ideal satisfying $(\mathrm{UR})$. The uniformity $\mathscr U_I$ has
a \emph{countable radius base} if there is a sequence $(r_n)_{n\geq1}$
in $I$ such that, for every $r\in I$, some $n$ satisfies
\begin{equation}\label{eq:6-3}
U_{r_n}\subseteq U_r.
\end{equation}
\end{definition}

\Needspace{12\baselineskip}
\begin{theorem}\label{thm:countability-new}
Let $L$ be a commutative unital quantale and let $I$ be a lattice
ideal satisfying $(\mathrm{UR})$. Then:
\begin{enumerate}
\item a countable radius base implies first countability $(A_1)$;
\item a countable radius base and $I$-total boundedness imply
second countability $(A_2)$;
\item a countable radius base and $T_0$ imply metrizability.
\end{enumerate}
In particular, compactness and a countable radius base imply
$(A_1)$ and $(A_2)$.
\end{theorem}

\begin{proof}
$(1)$ The balls $C_I(x,r_n)$ give a countable neighborhood base at $x$.

$(2)$ Let $(\rho_n)$ be a radius base and set $r_1=\rho_1$.
Choose $v_n\in I$ with $U_{v_n}\circ U_{v_n}\subseteq U_{r_n}$
and put $r_{n+1}=v_n\vee\rho_{n+1}$.
Since $r_n\geq\rho_n$, we have
$U_{r_n}\subseteq U_{\rho_n}$, so $(r_n)$ is a radius base.
Also, $r_{n+1}\geq v_n$ gives
$U_{r_{n+1}}\circ U_{r_{n+1}}
\subseteq U_{v_n}\circ U_{v_n}\subseteq U_{r_n}$. Hence
\begin{equation}\label{eq:6-4}
U_{r_{n+1}}\circ U_{r_{n+1}}\subseteq U_{r_n}\ (n\geq1).
\end{equation}
Choose finite $K_n$ with
$L=\bigcup_{a\in K_n}C_I(a,r_{n+1})$.
The family
\[
\mathscr B=\{C_I(a,r_{n+1}):n\geq1,\ a\in K_n\}
\]
is countable, since each $K_n$ is finite.
We show that it is a base for $\tau_I$.
Let $O$ be open and $x\in O$. Choose $n$ with
$C_I(x,r_n)\subseteq O$ and $a\in K_n$ such that
$x\in C_I(a,r_{n+1})$.
If $y\in C_I(a,r_{n+1})$, symmetry gives
$(x,y)\in U_{r_{n+1}}\circ U_{r_{n+1}}\subseteq U_{r_n}$
by \eqref{eq:6-4}. Therefore
\[
x\in C_I(a,r_{n+1})\subseteq C_I(x,r_n)\subseteq O,
\]
as required.

$(3)$ Under $T_0$, Theorem~\ref{thm:separation-new} gives
$\bigcap_{r\in I}U_r=\Delta_L$. A separated uniformity with a
countable base is metrizable \cite{Isbell1964,Kelley1975}.
The final assertion follows from Theorem~\ref{thm:compactness-new}.
\end{proof}

\section{Conclusion}

We have constructed ball topologies from lattice ideals on
commutative unital quantales without assuming involutivity.
Negation identifies ideal balls with filter balls in the
involutive integral case. Principal ideals on integral quantales
also admit multiplicative-filter representations, whereas
Theorem~\ref{ex:failure-new} shows that a general lattice ideal
need not have such a representation, even on a frame.

For integral quantales, continuity of multiplication at the unit
is equivalent to uniform ball refinement. Under these conditions,
all four operations are uniformly continuous.
Example~\ref{ex:power-new} shows that meet need not be continuous
without integrality. Example~\ref{ex:finite-ideal-JC-new} shows
that continuity of the finitary operations does not imply
join compatibility.

Two questions remain. First, which nonprincipal lattice ideals on
noninvolutive integral quantales admit a filter-ball representation?
Second, under uniform refinement, can completeness and total
boundedness be described by conditions on the order and residual
action on the ideal? For principal ideals, completeness is
automatic and compactness is equivalent to finiteness of the
quotient. A similar algebraic criterion for general ideals would
be useful.

\section*{Acknowledgements}

This work was supported by the National Natural Science Foundation
of China (Grant Nos. 12231007 and 12371462) and the Heilongjiang
Natural Science Foundation Joint Key Project (Grant No. ZL2024A001).

\begingroup
\small
\setlength{\bibsep}{0pt plus 2pt}

\endgroup

\end{document}